\documentstyle[xypic,epic,amsthm,amsfonts,amsbsy,amssymb,amsgen,amsmath,amsopn,verbatim,10pt]{article}
\input amssymb.sty
\input epsf

\newtheorem{theorem}{Theorem}[section]
\newtheorem{lemma}[theorem]{Lemma}
\newtheorem{proposition}[theorem]{Proposition}

\newtheorem{conjecture}[theorem]{Conjecture}
\newtheorem{definition}[theorem]{Definition}

\newtheorem{theorem-construction}[theorem]{Theorem - Construction}

\begin{document}

\newcommand{\Z}{{\Bbb Z}}
\newcommand{\R}{{\Bbb R}}
\newcommand{\Q}{{\Bbb Q}}
\newcommand{\C}{{\Bbb C}}
\newcommand{\lra}{\longrightarrow}
\newcommand{\lms}{\longmapsto}
\newcommand{\AAA}{{\Bbb A}}
\newcommand{\Alt}{{\rm Alt}}
\newcommand{\wg}{\wedge}
\newcommand{\ol}{\overline}
\newcommand{\CP}{{\Bbb C}P}
\newcommand{\bwg}{\bigwedge}
\newcommand{\caL}{{\cal L}}
\newcommand{\PP}{{\Bbb P}}
\newcommand{\HH}{{\Bbb H}}
\newcommand{\LL}{{\Bbb L}}

\begin{titlepage}
\title{Explicit regulator maps on polylogarithmic motivic complexes}
\author{A. B.  Goncharov}

\end{titlepage}
\date{}
\maketitle
\tableofcontents

\section  {Introduction}

{\bf 1. Regulator maps on the level of complexes}. Let $X$ be  an algebraic variety. Beilinson 
[B1] defined the rational motivic cohomology of $X$ 
via the algebraic K-theory of $X$ by the formula 
$$
H^i_{Mot}(X, \Q(n)):= gr^{\gamma}_nK_{2n-i}(X) \otimes \Q 
$$
Beilinson [B2]  and Lichtenbaum [L] conjectured that the weight $n$ 
motivic cohomology of $X$ should appear 
as cohomology groups of some complexes, called the weight $n$ motivic complexes of $X$.

If $X$ is defined over $\C$ Beilinson [B1] constructed the regulator map to the Deligne cohomology 
of $X$:
$$
H^i_{Mot}(X, \Q(n)) \lra H^i_{D}(X, \R(n))
$$
The regulator map plays key role 
in  the (hypothetical) formulas for special values of $L$-functions of 
motives over number fields [B1]. 

Our point is that the regulator map should be {\it explicitly} defined on the level of complexes. 
So for any algebraic variety over $\C$ one should have homomorphisms of complexes
\begin{equation} \label{3/1/00.1}
\mbox{weight $n$ motivic complex of $X$} \quad \lra \quad \mbox{weight $n$ Deligne complex of $X$}
\end{equation}
The cone of this map, shifted by $-1$, 
 defines the {\it Arakelov motivic complex}:
$$
{\rm R}\Gamma_X^{{\cal A}}(n):= {\rm Cone}\Bigl(\mbox{the map (\ref{3/1/00.1})}\Bigr)[-1]
$$
 and so its cohomology 
 should be called the {\it Arakelov motivic cohomology} of $X$. The group  
$H^{2n}({\rm R}\Gamma_X^{{\cal A}}(n))$ is canonically isomorphic to the group of codimension $n$ Arakelov cycles on $X(\C)$ (see s. 3.2). This isomorphism is transparent for a version of the Deligne complex 
recalled in s. 3.1. Examples are given in the 
section 3. (When $X$ is defined over $\Q$ we should take the cone of the regulator map to 
the so called real Deligne cohomology of $X \otimes_{\Q}\R$, see s. 3.2).

Motivic complexes are objects 
of the derived category. Several candidates for motivic complexes are known, 
each with its own charm. They  
 should be quasiisomorphic. Explicit regulator maps  should be 
defined for each of them. 

First of all there are Bloch's higher Chow complexes  ([Bl]). 
They satisfy many of the expected properties of motivic complexes. 
Explicit regulator maps from these complexes to the Deligne complexes were constructed 
in [G4] using the Chow polylogarithm construction given in 
[G3]. 

The goal of the present paper is to define regulator maps for another version  
of motivic complexes, the polylogarithmic complexes ([G1-2]). Our 
regulator maps  are defined very 
explicitly via the  classical polylogarithm functions
 with some funny combinations of Bernoulli numbers serving as the  coefficients. 
 
Combining this with Beilinson's conjecture on regulators we get, as a bonus,  
a precise conjecture on special values of $L$-functions of varieties over number 
fields. If the variety in question is spectrum of a number field it 
reduces precisely 
to Zagier's conjecture. So  our conjecture  is in the same relationship to 
Beilinson's conjecture as Zagier's conjecture [Z] to the 
Borel theorem [Bo]. 

The regulator maps for  the polylogarithmic complexes of weights $n \leq 3$ were 
constructed in [G1-2] and played an important role in the proofs of Zagier's and Deninger's 
conjectures on $\zeta$-functions at $s=3$ (loc. cit., [G5]).

{\bf 2. Polylogarithmic complexes}. Let $F$ be an arbitrary field. In [G1-2] we defined a complex $\Gamma(F;n)$ of the following shape
\begin{equation} \label{3/1/00.2}
{\cal B}_{n}\stackrel{\delta}{\rightarrow} {\cal B}_{n- 
1}\otimes F^{\ast} \stackrel{\delta}{\rightarrow} {\cal B}_{n- 
2}\otimes \Lambda^{2}F^{\ast}\stackrel{\delta}{\rightarrow} \ldots \stackrel{\delta}{\rightarrow} 
{\cal B}_{2}\otimes \Lambda^{n- 
2}F^{\ast}\stackrel{\delta}{\rightarrow} \Lambda^{n}F^{\ast} 
\end{equation}
called the polylogarithmic complex. Here the group ${\cal B}_{n}$ sits in degree $1$, 
and the differential is of degree $1$. We conjectured that this complex is quasiisomorphic to the weight $n$ motivic complex  of ${\rm Spec}(F)$, so one should have
\begin{equation} \label{3/1/00.3}
H^i(\Gamma(F;n))\otimes \Q \stackrel{?}{=} gr^{\gamma}_nK_{2n-i}(F) \otimes \Q 
\end{equation}
There is a good evidence this is so for small weights, see s. 1.3 below. 

Now let $F = \C(X)$ be the  field of rational functions on a complex algebraic variety $X$. 
Our main results  are explicit formulas 
for the regulator maps from the 
polylogarithmic complexes to the Deligne complexes of 
 ${\rm Spec}(\C(X))$, see theorems \ref{TCON} and  \ref{CANHH*}. The group ${\cal B}_n$ is directly related to the properties of the classical polylogarithms. In particular there is a homomorphism \begin{equation} \label{III}
{\cal B}_n(\C) \lra \R(n-1)
\end{equation}
 given in terms of the classical polylogarithms. 
Surprisingly its generalization to a homomorphism of complexes is  quite complicated. 

Our regulator maps enjoy 
{\it compatibility with the residues} property (condition (d) in theorem \ref{CANHH}), 
which would 
guarantee that they can be extended from the generic point of $X$ to $X$ itself. However there is a serious difficulty in the 
 definition of  polylogarithmic complex $\Gamma(X;n)$ 
for a general variety $X$ and $n>3$, (see p. 240 in [G1]). It would be resolved if homotopy invariance 
of complexes (\ref{3/1/00.2}) will be known (see conjecture 1.39 in [G1]). 
As a result we have unconditional definition of the polylogarithmic 
complexes $\Gamma(X;n)$ only in the following cases: 

a) $X = Spec(F)$, $F$ is an arbitrary field.

b) $X$ is an regular curve over any field, and $n$ is arbitrary.

c) $X$ is an arbitrary regular scheme, but $n \leq 3$. 

Compatibility with the residues provides 
the regulator map on the level of complexes in all these cases, 
and it would provide it in general if the mentioned above difficulty 
will be resolved. 

{\bf 3. Comparison with Beilinson's regulator map}. A homomorphism
from the motivic cohomology to $H^i(\Gamma(X;n))\otimes \Q $ 
has been constructed in the following cases:

(1) $F$ is an arbitrary field, $n \leq 3$ ([G1-2,5]) and $n=4$, $i>1$ (to appear).

(2) $X$ is a curve over a number field, $n \leq 3$ ([G5]) and $n=4$, $i>1$ (to appear). 

In all these cases we proved that this homomorphism followed by the regulator map 
on 
polylogarithmic complexes (when $F = \C(X)$ in (1)) coincides with Beilinson's regulator. 
This proves the traditionally difficult ``surjectivity'' property: 
the image of the regulator map on these polylogarithmic complexes {\it contains} 
the image of Beilinson's regulator map in the Deligne cohomology. 
 
The results (2) combined with the results of R. de Jeu [RdJ1-2] prove 
that the image of Beilinson's regulator map 
in the Deligne cohomology {\it coincides} with the image of the regulator map on polylogarithmic complexes in the case (2). 

{\bf 4. Remarks}. We give a detailed proof of the main result. It 
 is a rather involved but direct calculation. 
However the reader may skip it because   
there exists a  completely different, more conceptual,  approach to this result.  
We outline it in s. 2.8 to show the main theorem in a  general framework. 
The details 
will be published elsewhere. 
This approach needs quite elaborate and rather sofisticated machinery, so 
the elementary  proof 
presented in the last section might be the quickest   
way to check the main theorem.

{\bf Acknowledgment}.   
The author gratefully acknowledges the support
of the NSF grant  DMS-9800998. 
I am grateful to Spencer Bloch for discussions on Arakelov motivic complexes. 

I am delighted to contribute this paper 
to the volume telated to 20-th anniversary of Spencer Bloch's 
Irvine lectures, which always have been a  great source of inspiration for me.

\section  {The main result}

{\bf 1. Classical polylogarithms}. Recall their definition: 
$$
Li_1(z) := -\log(1-z); \qquad Li_n(z):=\int_0^z Li_{n-1}(t)
\frac{dt}{t}, \quad n\ge 2.
$$
They are multivalued analytic functions, but admit  the following 
single-valued cousins: 
$$
\widehat {\cal L}_n(z) := \quad \pi_n\Bigl(\sum_{k=0}^{n-1} \beta_k Li_{n-k}(z)\log^{n-k}|z|\Bigr)
$$
where
$$\pi_n(a + ib)=
\cases
a \qquad &\text{$n$ odd} \\
ib&\text{$n$ even}. 
\endcases$$
and
$$
\beta_{k}:=  \frac{2^kB_k}{k!}, \qquad \sum_{k \geq 0} \beta_{k}t^k = \frac{2t}{e^{2t} -1}
$$
so $\beta_{2m+1} = 0$ for $m\geq 1$, and  
$$
\beta_{0} = 1, \quad \beta_{1} = -1, \quad \beta_{2} = \frac{1}{3}, \quad \beta_{4} = -\frac{1}{45}, \quad \beta_{6} = \frac{2}{945}, ...
 $$
These functions were written by Zagier [Z]. Their Hodge-theoretic interpretation was given by Beilinson and Deligne [BD]. 
For example for the dilogarithm
it is the Bloch-Wigner function
$$ \widehat {\cal L}_2(z) := i{\cal L}_2(z):= 
\pi_2\left( Li_2(z)\right)+ i \arg(1-z)\cdot \log|z|$$

{\bf 2. The groups ${\cal B}_{n}(F)$ and 
polylogarithmic complexes(see s.1.4 in [G2])}.  For  a set $X$ denote by  $\Z[X]$ 
the free abelian group generated by symbols $\{x\}$ where $x$ run
through all elements of the set $X$. Let $F$ be an arbitrary field. 
We  define inductively subgroups
${\cal R}_{n}(F)$ of $ \Z[P^{1}_{F}]$, $n\geq 1$ and set 
$$
{\cal B}_{n}(F):=\Z[P^{1}_{F}]/{\cal R}_{n}(F) 
$$ 
By definition $ 
{\cal R}_{1}(F):=(\{x\} + \{y\} - \{xy\},(x,y\in F^{\ast}); 
\{0\};\{\infty\})\;$.  
Then ${\cal B}_{1}(F) = F^{\ast}$. 
 Let $\{x\}_{n}$ be the image of $\{x\}$ in ${\cal 
B}_{n}(F)$. Consider homomorphisms
\begin{eqnarray} 
 \mbox{$\Z$}[P^{1}_{F}] 
\stackrel{\delta_{n}}{\longrightarrow}  \left\{ 
\begin{array}{lll} 
{\cal B}_{n-1}(F)\otimes F^{\ast} &:& n\geq 3 \\ 
\Lambda^{2}F^{\ast} &:& n=2\end{array}\right. 
\end{eqnarray}
\begin{eqnarray}
 \delta_{n}:\{x\}\mapsto \left\{ \begin{array}{lll} 
\{x\}_{n-1}\otimes x &:& n\geq 3 \\ 
(1-x)\wedge x &:& n=2\end{array}\right.
 \qquad \delta_{n}:\{\infty\},\{0\},\{1\}\mapsto 0 
\end{eqnarray} 
  Set 
$ {\cal A}_{n}(F):={\rm Ker}\ \delta_{n}\; . $ 
 Any element $\alpha(t) = \Sigma n_{i}\{f_{i}(t)\} \in \Z[P^{1}_{F(t)}]$ has a specialization $\alpha(t_{0}):=\Sigma 
n_{i}\{f_{i}(t_{0})\}\in \Z[P^{1}_{F}]$ at each point $t_{0}\in 
P^{1}_{F}$.  
 
\begin {definition}
\label {2.10}  
${\cal R}_{n}(F)$  is generated by 
elements $\{\infty\},\{0\}$  and   $\alpha(0)-\alpha(1)$  where $\alpha(t)$ runs 
through all elements of ${\cal A}_{n}(F(t))$.
\end {definition}
 Then $\delta_{n}\Bigl({\cal R}_{n}(F)\Bigl)=0$ ([G1], 1.16).  So we get 
homomorphisms  
$$
\delta_n: {\cal B}_{n}(F) \longrightarrow 
{\cal B}_{n-1}(F)\otimes F^{\ast}, \quad n\geq 3; \quad  \delta_2: {\cal
B}_{2}(F) \longrightarrow  \Lambda^{2}F^{\ast} 
$$
and finally the mentioned in the introduction complex $\Gamma (F,n)$: 
$$
{\cal B}_{n}\stackrel{\delta}{\rightarrow} {\cal B}_{n- 
1}\otimes F^{\ast} \stackrel{\delta}{\rightarrow} {\cal B}_{n- 
2}\otimes \Lambda^{2}F^{\ast}\stackrel{\delta}{\rightarrow} \ldots \stackrel{\delta}{\rightarrow} 
{\cal B}_{2}\otimes \Lambda^{n- 
2}F^{\ast}\stackrel{\delta}{\rightarrow} \Lambda^{n}F^{\ast} 
$$ 
where   
$\delta:\{x\}_{p}\otimes \bigwedge^{n-p}_{i=1} y_{i}\to 
\delta_p(\{x\}_{p})\wedge \bigwedge^{n-p}_{i=1}y_{i} $.

Let $F = \C$. 
Set $\widehat {\cal L}_{n}(\sum m_i \{z_i\}_n):= \sum m_i \widehat {\cal L}_{n}(z_i)$. 
One can prove that $\widehat {\cal L}_{n}\Bigl({\cal R}_{n}(\C))\Bigr)=0$ ([G2], theorem 1.13). 
So we are getting homomorphism (\ref{III}). 

{\bf 3. The residue homomorphism for complexes $\Gamma(F,n)$ (s. 1.14 in [G1])}.
Let $F=K$ be a  field with 
a discrete valuation $v$,  the residue field  $k_v$ and  the group of units
$U$. Let $u \rightarrow \bar u$ be the projection $U \rightarrow
k_v^{\ast}$. Choose  a uniformizer $\pi$. There is a homomorphism $\theta:\Lambda^{n}K^{\ast} 
\longrightarrow\Lambda^{n-1} k_v^{\ast}$ uniquely defined 
by the  following properties $(u_{i}\in U)$: 
$$
\theta\; (\pi\wedge u_{1}\wedge \cdots\wedge u_{n-1}) 
= \bar u_{1}\wedge\cdots \wedge \bar u_{n-1}; \qquad  \theta\;
(u_{1}\wedge \cdots \wedge u_{n}) = 0 
$$
It is clearly independent  of $\pi$. Define a homomorphism $s_{v}:\Z[P^{1}_{K}]\longrightarrow \Z[P^{1}_{ k_v}]$ by setting
 $s_{v}\{ x\} =  \{ \bar x\}  \mbox{ if $x$ is a unit}$ and $0$  
 otherwise.
It induces a homomorphism 
$s_{v} :  {\cal B}_{m}(K)\longrightarrow {\cal B}_{m}(k_v)$. 
Put
$$ 
\partial_{v}:= s_{v}\otimes \theta : \quad {\cal B}_{m}(K)\otimes 
\Lambda^{n-m} K^{\ast}\quad \longrightarrow \quad {\cal B}_{m} (k_v) 
\otimes \Lambda^{n-m-1} k_v^{\ast} 
$$ 
It  defines a morphism 
of complexes
$\partial_{v} : \Gamma (K,n)\longrightarrow  \Gamma( k_v,
n-1)[-1]$.

{\bf 4. Main result: a preliminary form}. Let ${\cal A}^i(\eta_X)$ be the space of real 
smooth $i$-forms at the generic point $\eta_X:= {\rm Spec} \C(X)$ of a 
complex variety $X$. Denote by 
  $X^{(1)}$ the set of the codimension one closed irreducible subvarieties in $X$. Let 
${\cal D}$ be  
the de Rham differential on distributions, and $d$ the de Rham differential on 
${\cal A}^i(\eta_X)$. A typical example:
\begin{equation} \label{3/5/00.6}
d \Bigl( d i \arg z\Bigl) = 0; \quad {\cal D}\Bigl( d i \arg z\Bigl) = 2 \pi i \delta(z)
\end{equation} 
The difference ${\cal D} - d$ is the de Rham residue homomorphism. It is defined 
on distributions 
smooth at the generic point of $X(\C)$. Its value on such a distribution is 
concentrated on a union of codimension one subvarieties. 
 
\begin{theorem} \label{CANHH}
\label {1.11a} There exist a homomorphism of complexes
$$
\begin{array}{ccccccc}
{\cal B}_{n}(\C(X))&\stackrel{\delta}{\rightarrow}& {\cal B}_{n- 
1}(\C(X))\otimes \C(X)^{\ast} & \stackrel{\delta}{\rightarrow} & \ldots &\stackrel{\delta}{\rightarrow} &
 \bigwedge^{n}\C(X)^{\ast}\\ 
&&&&&&\\
\downarrow r_n(1)& &\downarrow r_n(2)& & & &\downarrow r_n(n)\\
&&&&&&\\
{\cal A}^0(\eta_X)(n-1) &\stackrel{d}{\rightarrow}& {\cal A}^1(\eta_X)(n-1) &\stackrel{d}{\rightarrow}& ...
&\stackrel{d}{\rightarrow}& {\cal A}^{n-1}(\eta_X)(n-1)\\
\end{array}
$$
such that 

(a) $r_n(1)(\{f\}_n) = \widehat{{\cal L}}_n(f)$.

(b) $ d r_n(n)(f_1 \wedge ... \wedge f_n) + \pi_n (d\log f_1 \wedge
... \wedge d\log f_n)
= 0$.

(c) The differential form $r_n(m)(\ast)$ defines a distribution on $X(\C)$. 

(d) The homomorphism $r_n(\cdot)$ is compatible with residues:
$$
{\cal D} \circ r_n(m) - r_n(m+1) \circ \delta = 2 \pi i \cdot \sum_{Y \in X^{(1)}}r_{n-1}(m-1) \circ \partial_{v_Y}, \quad m<n
$$
$$
{\cal D} \circ r_n(n) - \pi_n(d\log f_1 \wedge
... \wedge d\log f_n) = 
2 \pi i\cdot \sum_{Y \subset X^{(1)}}r_{n-1}(n-1) \circ \partial_{v_Y}
$$
where $v_Y$ is  the  valuation  
on the field $\C(X)$ defined by a divisor $Y$.
 \end {theorem}

 {\bf Remark}. This result has been formulated in [G4], see theorem 4.3.

The part (d) means that  $r_n(\cdot)$  sends the   residue
homomorphism $\partial_{v_Y}$
 to the de Rham  residue homomorphism.

Here are two  examples. Set
$$
\alpha(f,g):= -\log|f|d\log|g| + \log|g|d\log|f|
$$

{\bf Example 1: n=3}.
Set
\begin{eqnarray*}
& & r_{3}(2) : \{ f\}_{2} \otimes g\longmapsto \widehat{{\cal L}}_{2} (f) d i\arg g -\frac{1}{3} \log \vert g\vert \cdot
\alpha(1-f,f)\\  
\end{eqnarray*}

{\bf Example 2: n=4}.
\begin{eqnarray*}
 & &r_{4}(2) : \{ f \}_{3} \otimes g \mapsto  \widehat{{\cal L}}_{3} (f ) d i \arg g    -\frac{1}{3} \widehat{{\cal L}}_{2} (f )\log \vert g \vert \cdot d \log \vert f\vert\\
& & r_{4}(3) : \{ f \}_{2} \otimes g_1 \wedge g_2 \mapsto + \widehat{{\cal L}}_{2} (f ) d i\arg g_1 \wedge d i\arg g_2   -
\frac{1}{3} \alpha(1-f,f) \cdot \\
& &\Bigl(\log \vert g_1 \vert  d i\arg \vert g_2\vert  -
\log \vert g_2 \vert  d i\arg \vert g_1\vert\Bigr) \quad 
+ \frac{1}{3} \widehat{{\cal L}}_{2} (f ) d \log \vert g_1 \vert \wedge d \log \vert g_2 \vert\\
\end{eqnarray*}

{\bf 5. The numbers $\beta_{k,p}$}. 
Define for any integers $p \geq 1$ and $k \geq 0$ the numbers 
$$
\beta_{k, p}:= \quad (-1)^p (p-1)! \sum_{0 \leq i \leq [\frac{p-1}{2}]} \frac{1}{(2i+1)!} \beta_{k+p-2i}
$$
For instance
$$
\beta_{k, 1} = -\beta_{k+ 1}; \quad \quad \beta_{k, 2} = \beta_{k+ 2}; \quad \quad \beta_{k, 3} = 
-2(\frac{1}{1!}\beta_{k+ 3} + \frac{1}{3!}\beta_{k+ 1}); 
$$
$$
\beta_{k, 4} = 3!(\frac{1}{1!}\beta_{k+ 4} + \frac{1}{3!}\beta_{k+ 2}); \quad \quad 
\beta_{k, 5} = -4!(\frac{1}{1!}\beta_{k+ 5} + \frac{1}{3!}\beta_{k+ 3} + \frac{1}{5!}\beta_{k+ 1})
$$
One has  recursions
\begin{equation} \label{12.99}
2p\cdot \beta_{k+1, 2p} =  - \beta_{k, 2p+1} -\frac{1}{2p+1}\beta_{k+1}; \qquad 
(2p - 1)\cdot \beta_{k+1, 2p-1} =  - \beta_{k, 2p} 
\end{equation}
These recursions together with  $\beta_{k, 1} = -\beta_{k+ 1}$ determine the numbers $\beta_{k, p}$.

\begin{lemma} \label{RADEMA} Let $m \geq 1$. Then 
\begin{equation} 
\beta_{0, 2m} = \beta_{0, 2m+1} =  \frac{1}{2m+1} 
\end{equation}
\begin{equation} \label{1999.5}
\beta_{1, 2m-1}  =  -\frac{1}{(2m-1)(2m+1)}, \qquad \beta_{1, 2m}  = 0
\end{equation}
\end{lemma}

{\bf Proof}. Let us prove formula
\begin{equation} \label{FFOO}
\beta_{0, 2m} :=(2m-1)! \Bigl(\frac{1}{1!}\beta_{2m} +  \frac{1}{3!}\beta_{2m- 2} + ... + \frac{1}{(2m-1)!}\beta_{2}\Bigr) =
 \frac{1}{2m+1}
\end{equation}
Since $\beta_0 =1$, $\beta_1 =-1$ and $\beta_{2m+1} = 0$ for $m>0$ this formula, as one easily checks, 
is equivalent to the following one:
$$
\sum_{p \geq 1}\frac{1}{p!} \beta_{2m+1-p} = 0
$$
The left hand side is the coefficient in $t^{2m+1}$ of the power series 
$$
(e^t-1)\frac{2t}{e^{2t}-1} = \frac{2t}{e^{t} + 1} 
$$
The right hand side of the last equation is an almost even function in $t$: if we denote it by $F(t)$ then $F(-t) = F(t) + 2t$. 

Formula (\ref{FFOO}) immediately implies that 
$$
\beta_{1, 2m-1} := - (2m-2)! \Bigl(\frac{1}{1!}\beta_{2m} +  \frac{1}{3!}\beta_{2m- 2} + ... + \frac{1}{(2m-1)!}\beta_{2}\Bigr) =
 \frac{1}{(2m-1)(2m+1)}
$$

Using $\beta_{2m+1}=0$ for $m>0$ and $\beta_{1}= -1$ we have 
$$
\beta_{0, 2m+1} =-(2m)! \Bigl(\frac{1}{1!}\beta_{2m+1} +  \frac{1}{3!}\beta_{2m- 1} + ... + \frac{1}{(2m-1)!}\beta_{1}+ \frac{1}{(2m+1)!}\beta_{1}\Bigr) =
 \frac{1}{2m+1}
$$
and
$$
\beta_{1, 2m} =(2m-1)! \Bigl(\frac{1}{1!}\beta_{2m+1} +  \frac{1}{3!}\beta_{2m- 1} + ... + \frac{1}{(2m-1)!}\beta_{3}\Bigr) = 0
$$
The lemma is proved.

\begin{proposition} \label{3/5/00.3}For any $p \geq 1$ one has 
\begin{equation} \label{3/5/00.1}
\beta_{n-2,p+1} - (n-1)\beta_{n-1,p} - \sum_{k=1}^{n-3}\beta_{k,p} \beta_{n-k-1}  = 0
\end{equation}
\end{proposition}

{\bf Proof}. We will do it by induction on $p$ using the recursion relations (\ref{12.99}). 
If $p=1$ the formula we need to prove boils down to the identity
\begin{equation} \label{3/5/00.2}
\sum_{k=2}^{n-2}\beta_{i}\beta_{n-i} +n \beta_n =0
\end{equation}
which is easy to check using the generating functions. 
Let us denote by $(*)_{n,p}$ the left hand side of (\ref{3/5/00.1}). 
Then $(*)_{n,p} +(*)_{n+1,p-1} =0$. Indeed, 
one has  
\begin{equation} \label{3/5/00.2}
(*)_{n,p} +(*)_{n+1,p-1} \quad  = \quad \beta_{n-2,p+1} + (p-1)\beta_{n-1,p} 
\end{equation}
$$
-(n-1)\beta_{n-1,p} - n (p-1)\beta_{n,p-1} \quad - \quad \sum_{k=1}^{n-3}
\beta_{k,p}\beta_{n-k-1} - (p-1)\sum_{k=0}^{n-3}
\beta_{k+1,p-1}\beta_{n-k-1}  
$$
Let us assume first that $p$ is even. Using  recursions (\ref{12.99})
 we write (\ref{3/5/00.2}) as 
$$
-\frac{1}{p+1}\beta_{n-1} - \beta_{n-1,p} - (p-1)\beta_{n,p-1} - (p-1)\beta_{1,p-1}
\beta_{n-1}  
$$
To prove that it is zero we use  again recursions as well as the first formula 
in (\ref{1999.5}). 
Now let $p$ be an odd number. Then the recursion relations show that (\ref{3/5/00.2}) equals to 
$$
-\beta_{n-1,p} - (p-1)\beta_{n,p-1} +\frac{n-1}{p} \beta_{n} + \frac{1}{p}
\sum_{k=1}^{n-3}\beta_{k+1}\beta_{n-k-1}  - (p-1)\beta_{1,p-1}\beta_{n-1}
$$
Using the identity (\ref{3/5/00.2}) together with the recursions and the second formula 
in (\ref{1999.5}) we see that this expression is also equal to zero. 
The proposition is proved.

{\bf 6. Construction of the homomorphism $r_n(\cdot)$.}  
Let us define differential $1$-forms $\widehat {\cal L}_{p,q}$ on $\CP^1 \backslash \{0, 1, \infty\}$ for $q\geq 1$  as follows:
\begin{equation} \label{TG10}
\widehat {\cal L}_{p,q}(z):= \widehat {\cal L}_{p}(z) \log^{q-1}|z| \cdot d\log |z|, \quad p \geq 2
\end{equation}
$$
\widehat {\cal L}_{1,q}(z):=  \alpha(1-z, z) \log^{q-1}|z|
$$
It provides a distribution on $\CP^1$. Moreover for any rational function $f$ on a complex variety $X$ 
the $1$-form $\widehat {\cal L}_{p,q}(f)$ provides  a distribution on $X(\C)$.

{\it A useful notation}. Set 
$$
{\cal A}_m\left\{\bigwedge_{i=1}^{2p}d\log|g_i| \wedge \bigwedge_{i=2p+1}^{m}di \arg g_{j}
  \right\}:= 
$$
$$
{\rm Alt}_m\left\{\frac{1}{(2p)!(m-2p)!}\bigwedge_{i=1}^{2p}d\log|g_i| \wedge \bigwedge_{i=2p+1}^{m}di \arg g_{j} \right\} 
$$
and 
$$
{\cal A}_m\left\{\log|g_1| \cdot \bigwedge_{i=2}^{p}d\log|g_i| \wedge \bigwedge_{i=p+1}^{m}di \arg g_{j} \right\}:= 
$$
$$
{\rm Alt}_m\left\{\frac{1}{(p-1)!(m-p)!}\log|g_1|  \cdot \bigwedge_{i=2}^{p}d\log|g_i| \wedge \bigwedge_{i=p+1}^{m}di \arg g_{j}   \right\}
$$
So
${\cal A}_m(F(g_1, .., g_m))$ is a weighted alternation 
(we divide by the order of the stabilizer of the term we alternate). 

Now we are ready for the precise formulation of our main result. 

\begin{theorem-construction} \label{TCON}
Let $f, g_1, ..., g_m$ be rational functions on a complex variety $X$. 
Then the  following formula provides maps satisfying all the properties of theorem \ref{CANHH}. 
$$
r_{n+m}(m+1):\{f\}_n \otimes g_1 \wedge ... \wedge g_m \lms 
$$
\begin{equation} \label{12.26.11}
\widehat {\cal L}_{n}(f) \cdot {\cal A}_m\left\{  \sum_{p \geq 0 } \frac{1}{2p+1}
 \bigwedge_{i=1}^{2p} d\log|g_i| \wedge \bigwedge_{j=2p+1}^{m} d i \arg g_{j}
 \right\} +
\end{equation}
\begin{equation} \label{12.26.12}
\sum_{k \geq 1 }\sum_{1 \leq p \leq m} \beta_{k, p} \widehat {\cal L}_{n-k,k}(f) \wedge 
{\cal A}_m\left\{\log|g_1| \bigwedge_{i=2}^p d \log|g_{i}| \wedge \bigwedge_{j= p+1}^m 
d i \arg g_{j}
 \right\} 
\end{equation}
\end{theorem-construction}

Here are several examples.

{\bf Example 1}. $m=1$, $n$ is arbitrary.
$$
\{f\}_n \otimes g \lms \quad \widehat {\cal L}_{n}(f) d i \arg g
- \sum_{k=1}^{n-1} \beta_{k+1} \widehat {\cal L}_{n-k, k}(f) \cdot \log|g| 
$$

{\bf Example 2}. $m=2$, $n$ is arbitrary. 
$$
\{f\}_n \otimes g_1 \wedge g_2 \lms \widehat {\cal L}_{n}(f) \left\{  d i \arg g_1 \wedge d i \arg g_2 + 
\frac{1}{3} d \log |g_1| \wedge d \log |g_2|
\right\}
$$
$$
- \sum_{k=1}^{n-1} \beta_{k+1} \widehat {\cal L}_{n-k, k}(f) \wedge (\log|g_1| d i\arg g_2 -\log|g_2| d i\arg g_1)  
$$
$$
+ \sum_{k \geq 1} \beta_{k+2} \widehat {\cal L}_{n-k, k}(f) \wedge (\log|g_1| d \log|g_2| -\log|g_2| d \log|g_1|) 
$$

{\bf Example 3}. The homomorphism $r_5(*)$. 
$$
\{f\}_4 \otimes g \lms \quad \widehat {\cal L}_{4}(f) d i \arg g - \frac{1}{3} \widehat {\cal L}_{3}(f) 
d \log |f| \cdot \log|g| + \frac{1}{45} \alpha(1-f,f) \log^2|f| \cdot \log|g|
$$

$$
\{f\}_3 \otimes g_1 \wedge g_2 \lms \quad \widehat {\cal L}_{3}(f) \left\{d i \arg g_1 \wedge d i \arg g_2 +
 \frac{1}{3} d \log |g_1| \wedge d\log|g_2| \right\}
$$
$$
- \frac{1}{3} \widehat {\cal L}_{2}(f) d \log|f| \wedge \Bigl(\log|g_1| d i \arg g_2 - \log|g_2| d i \arg g_1\Bigr)
$$
$$
- \frac{1}{45} \alpha(1-f,f) \log |f| \wedge  \Bigl(\log |g_1| d \log |g_2| - \log|g_2| d \log |g_1|\Bigr)
$$ 

$$
\{f\}_2 \otimes g_1 \wedge g_2 \wedge g_3 \lms \quad 
$$
$$
\widehat {\cal L}_{2}(f) \cdot 
{\cal A}_3\Bigl(d i \arg g_1 \wedge d i \arg g_2 \wedge d i \arg g_3 +
\frac{1}{3} d \log |g_1| \wedge d\log|g_2| d i \arg g_3 \Bigr)
$$
$$
- \alpha(1-f,f) \wedge  {\cal A}_3\Bigl(\frac{1}{3} \log|g_1| d i \arg g_2 \wedge d i \arg g_3 
+ \frac{1}{15}  \log|g_1| d \log |g_2| \wedge d \log |g_3|\Bigr)
$$

{\bf Remark}. The morphism of complexes $r_n(\cdot)$ is 
not defined uniquely by its properties if $n>3$. For example if $n=4$ we can have a map homotopic to our regulator map by using the homotopy which is given by the homomorphism
$$
{\cal B}_2(\C(X)) \otimes \Lambda^2\C(X)^* \lra {\cal A}^0(\eta_X), \qquad \{f\}_3 \otimes g \lms 
\widehat {\cal L}_{2}(f)\cdot \log|f| \log|g|
$$
(the other components of the homotopy are zero). Indeed, the residue map for it is zero, and it takes values in $\R(2)$.

{\bf 7. Theorem \ref{TCON} from the point of view of the Deligne complex}.  
Recall that a  $p$-distribution on a manifold $X$ is 
 a linear continuous  functional on the space of 
$(\dim_{\R}X-p)$-forms with compact support.  Denote by 
${\cal D}^{p}_{X}$ the space of all real $p$-distributions on $X$.

Let $X$ be a regular variety over $\C$.  The $n$-th Beilinson-Deligne complex 
$\underline {\R}_{{\cal D}}(n)_X$ can be 
defined as a total complex associated with the following
bicomplex of sheaves in classical topology on $X(\C)$: 
$$
\begin{array}{ccccccccccc} \label{del}
\Bigl({\cal D}_{X}^{0}&\stackrel{d}{\longrightarrow}&{\cal
D}_{X}^{1}&\stackrel{d}{\longrightarrow}&\ldots&\stackrel{d}{\longrightarrow}&{\cal
D}^{n}_{X}&\stackrel{d}{\longrightarrow}&{\cal
D}_{X}^{n+1}&\stackrel{d}{\longrightarrow}&\ldots\Bigr) \otimes \R(n-1)\\
&&&&&&&&&&\\
&&&&&&\uparrow\pi_{n} &&\uparrow\pi_{n}&&\\
&&&&&&&&&&\\
&&&&&&\Omega^{n}_{X, \log}
&\stackrel{\partial}{\longrightarrow}&\Omega_{X, \log}^{n+1}&\stackrel{\partial}{\longrightarrow}&
\end{array}
$$
Here ${\cal D}^{0}_{X}$ placed in
degree 1 and  
$(\Omega^{\bullet}_{X, \log}, \partial)$ is the de Rham complex of 
holomorphic forms with logarithmic singularities
 at infinity. We will denote by $\underline 
{\R}_{{\cal D}}(n)(U)$ the complex of the global sections. 

Theorem \ref{CANHH} can be reformulated as follows

\begin{theorem} \label{CANHH*}
Let $X$ be a complex algebraic variety. Set $\widetilde r_n(i):= r_n(i)$ for $i<n$ and let 
$$
\widetilde r_n(n): \Lambda^n\C(X)^* \lra {\cal A}^{n-1}(\eta_X)(n-1) \oplus 
\Omega^n_{\log}(\eta_X)
$$
\begin{equation} \label{3/6/00.01}
f_1 \wedge ... \wedge f_n \lms r_n(n)(f_1 \wedge ... \wedge f_n) + d\log f_1 \wedge ... \wedge 
d\log f_n
\end{equation}

Then we get a homomorphism of complexes
\begin{equation} \label{CANHH**}
\widetilde r_n(\cdot): \Gamma(\C(X); n) \lra \underline {\R}_{{\cal D}}(n)(\eta_X)
\end{equation}
compatible with the residues 
as explained in the part (d) of theorem \ref{CANHH}. 
 \end {theorem}

Indeed, condition b) of theorem \ref{CANHH} just means that the right hand side of 
(\ref{3/6/00.01}) is a cycle 
in the Deligne complex $\underline {\R}(n)_{{\cal D}}(\eta_X)$. 

There is a natural Dolbeault resolution of the complex $\Omega_{\log}^{\geq n}$. Using it 
in the complex of sheaves $\underline {\R}_{{\cal D}}(n)_X$ to replace the subcomplex 
$\Omega_{\log}^{\geq n}$ we get a complex $\underline {{\R}}_{{\cal D}}'(n)_X$ 
of fine sheafs on $X(\C)$. The property (d) would allow us to extend the homomorphism 
$\widetilde r_n(\cdot)$ to a morphism of complexes
$$
\widetilde r_n(\cdot): \Gamma(X; n) \lra \underline { {\R}}'_{{\cal D}}(n)(X(\C))
$$
For small weights this was explained in detail in [G2]. 

If $X$ is a variety over $\R$, then
$$
H_{{\cal D}}^i(X_{\R}; \R(n)) = H_{{\cal D}}^i (X(\C), \R( n))^{\bar F_{{\infty}}}
$$
where $\bar F_{{\infty}}$ is the de Rham involution,  
i.e.  the composition  of the involution $F_{{\infty}}$   
on $X(\C)$ induced by the complex conjugation  
with the complex conjugation of coefficients. 
If $X$ is a variety over $\Q$ then the regulator map is defined as 
\begin{equation} \label{RRGG}
H^i\Gamma(X;n) \lra H^i_{{\cal D}}(X\otimes_{\Q} {\R}/_{\R}; \R(n))
\end{equation}

\begin{conjecture} \label{MMMCCC}
The image of the regulator map (\ref{RRGG}) in the Deligne cohomology 
coincides with the image of Beilinson's regulator map.  
\end{conjecture}

If $X = {\rm Spec}(F)$ where $F$ is a number field then the only nontrivial case is $i=1$, and we 
get a version of Zagier's conjecture. 
The next case is when $X$ is a curve over number field, and $i=2$. Then we come to  
conjecture 1.5 in [G5], see also theorem 4.4 in [G5] and 
s. 7 in [G7]. 

There is a  version of conjecture \ref{MMMCCC} expressing the special value of 
the corresponding $L$-function at $s=n$ via the regulator map on the polylogarithmic complex. 
Its specialization for elliptic curves is conjecture 1.10 in [G7], see also theorem 4.7 in [G5]. 
For the $n=2$ it is a theorem of Bloch for elliptic CM curves and of Beilinson 
for elliptic curves over $\Q$, for $n=3$ it was conjectured by Deninger and proved 
by the author in [G5] for elliptic curves over $\Q$, and it is also proved now by the author 
for $n=4$ for 
elliptic curves over $\Q$ (to appear).

{\bf 8. Generalizations}. According to the Tannakian formalism the category 
of mixed $\R$-Hodge-Tate structures is equivalent to the category of graded comodules over a certain  Lie coalgebra 
${\cal L}_{\bullet}(\C) $ positively graded by the weights. 
There exists a natural homomorphism of groups 
$$
{\cal B}_{n}(\C) \lra {\cal L}_{n}(\C)  
$$
provided by the Hodge realization of the polylogarithmic motive. One has a 
homomorphism of $\Q$-vector spaces 
$$
 p_n: {\cal L}_{n}(\C)  \lra \R(n-1)
$$
called the Lie-period map, see [D] for its construction. The  composition 
$$
{\cal B}_{n}(\C) \lra {\cal L}_{n}(\C)  \lra \R(n-1)
$$
coincides with the homomorphism $\widehat {\cal L}_n$, see [BD]. 

We generalize this picture considering variations of mixed Hodge-Tate structures over complex varieties and working on the level of  complexes. 
 
First of all,  the category 
of variations of mixed $\R$-Hodge-Tate structures over $\eta_X= {\rm Spec}\C(X)$, where $X$ is a complex algebraic 
variety,  
is equivalent to the category of graded comodules over a  Lie coalgebra 
${\cal L}_{\bullet} (\eta_X)$.

The Hodge realization functor provides canonical homomorphism 
\begin{equation} \label{3/7/00.1}
l_n: {\cal B}_{n}(\C(X)) \lra {\cal L}_{n} (\eta_X)
\end{equation}
Applying it componentwise  we get a morphism of complexes ([G1-2]) 
$$
l_{n-k} \otimes \wedge^k l_1: {\cal B}_{n- k}(\C(X))\otimes \Lambda^k\C(X)^{\ast} \lra 
{\cal L}_{n-k} (\eta_X) \otimes \Lambda^k{\cal L}_{1} (\eta_X) 
$$
The vector space on the right is a subspace in 
$
 \Bigl( \Lambda^{k+1}{\cal L}_{\bullet} (\eta_X)\Bigr)_n
$ 
where the subscript $n$ means the weight $n$ part. 
The standard cochain complex of a Lie coalgebra 
${\cal L}_{\bullet} $ looks as follows
\begin{equation} \label{6/3/00.2}
{\cal L}_{\bullet} (\eta_X) \lra \Lambda^2{\cal L}_{\bullet} (\eta_X) \lra 
 \Lambda^3{\cal L}_{\bullet} (\eta_X) \lra ...  
\end{equation}
It is a complex of graded $\Q$-vector spaces. We get 
a morphism of complexes 
$$
\begin{array}{ccccccc}
{\cal B}_{n}(\C(X))&\stackrel{\delta}{\rightarrow}& {\cal B}_{n- 
1}(\C(X))\otimes \C(X)^{\ast} & \stackrel{\delta}{\rightarrow} & \ldots &\stackrel{\delta}{\rightarrow} &
 \bigwedge^{n}\C(X)^{\ast}\\ 
&&&&&&\\
\downarrow & &\downarrow & & & &\downarrow \\
&&&&&&\\
{\cal L}_{n} (\eta_X) &\stackrel{}{\rightarrow}& 
\Bigl(\Lambda^2{\cal L}_{\bullet} (\eta_X)\Bigr)_n &
\stackrel{}{\rightarrow}& ...
&\stackrel{}{\rightarrow}& \Lambda^n{\cal L}_{1} (\eta_X)\\
\end{array}
$$
The bottom line of this diagram is the degree $n$ part of the complex (\ref{6/3/00.2}). 

The Lie-period provides a natural homomorphism of $\Q$-vector spaces 
\begin{equation} \label{6/3/00.19}
 p_n: {\cal L}_{n} (\eta_X) \lra {\cal A}^0(\eta_X)(n-1)
\end{equation}
The following result tells us
 that it is a beginning of a  morphism of complexes. 

\begin{theorem} \label{6/3/00.1} There exists a morphism of  complexes
\begin{equation} \label{6/3/00.20}
\mbox{the weight $n$ part of
 (\ref{6/3/00.2})} \lra \mbox{the weight $n$ Deligne complex 
$\underline {\R}_{{\cal D}}(n)(\eta_X)$}
\end{equation}
 whose degree $1$ component
$ {\cal L}_{\bullet}(\eta_X) \lra {\cal A}^0(\eta_X)$ coincides 
with the map (\ref{6/3/00.19}). 
\end{theorem}
A detailed account on this result and related issues would double the size of this paper.
We will do it in a different place. 

Combining homomorphism  (\ref{6/3/00.20}) with the previous map of complexes 
$$
\Gamma(\eta_X;n) \lra \mbox{the weight $n$ part of 
 (\ref{6/3/00.2})} 
$$
we get 
a homomorphism of complexes from theorem \ref{CANHH*}.

\section  {Arakelov motivic complexes: examples}

{\bf 1.  Another version of  Deligne's complexes}.
We need a complex   introduced twenty years ago by Deligne and  
 quasiisomorphic to the complex $\underline {\R}'_{{\cal D}}(n)(X(\C))$.

Let ${\cal D}_X^{ p,q} = {\cal D}^{ p,q}$ be the abelian group of complex valued distributions of type $(p,q)$ on $X(\C)$. Consider the following cohomological 
bicomplex,  where the group ${\cal D}^{0,0}$ is placed in degree $1$:
$$
\begin{array}{ccccccccc}
&&&&&&&&{\cal D}_{cl}^{n,n}\\
&&&&&&&&\\
&&&&&&&2 \bar \partial \partial\nearrow&\\
&&&&&&&&\\
{\cal D}^{0,n-1}&\stackrel{\partial}{\longrightarrow}&{\cal D}^{1,n-1}&\stackrel{\partial}{\longrightarrow}&
...&\stackrel{\partial}{\longrightarrow}&{\cal D}^{n-1,n-1}&&\\
&&&&&&&&\\
\bar \partial \uparrow &&\bar \partial \uparrow &&&&\bar \partial \uparrow&&\\
 ...&...&...&...&...&...&...&&\\
 \bar \partial \uparrow&&\bar \partial \uparrow&&&&\bar \partial \uparrow&&\\
&&&&&&&&\\
{\cal D}^{0,1}&\stackrel{\partial}{\longrightarrow}& {\cal D}^{1,1}&\stackrel{\partial}{\longrightarrow}&... &\stackrel{\partial}{\longrightarrow}& {\cal D}^{n-1,1}&&\\
&&&&&&&&\\
\bar \partial \uparrow&&\bar \partial \uparrow&&&&\bar \partial \uparrow&&\\
&&&&&&&&\\
{\cal D}^{0,0}&\stackrel{\partial}{\longrightarrow} &{\cal D}^{1,0}& \stackrel{\partial}{\longrightarrow} &...&\stackrel{\partial}{\longrightarrow} &{\cal D}^{n-1,0}&&
\end{array}
$$
 
Let $Tot^{\bullet}$ be the total complex of this  bicomplex. It is concentrated in degrees $[1,2n]$. The   complex $C^{\bullet}_{{\cal D}}(X; \R(n))$ is  subcomplex of $Tot^{\bullet}$ 
defined as follows. 
Intersect the part of the complex $Tot^{\bullet}$  
  coming from the $n \times n$ square in the diagram with the complex of $\R(n-1)$-valued 
distributions. 
Consider the   subgroup    ${\cal D}_{cl, \R}^{n,n}(n) \subset {\cal D}_{cl}^{n,n}$  of
 the $\R(n)$-valued distributions of type $(n,n)$. They  form a subcomplex in $Tot^{\bullet}$ because
 $\bar \partial \partial$ sends $\R(n-1)$-valued distributions to $\R(n)$-valued distributions.
This is the complex $C^{\bullet}_{{\cal D}}(X; \R(n))$ with the differential denoted by $D$. 
It is a truncation of the complex   considered by Deligne.

\begin{proposition} \label{dcom}
The   complex $C^{\bullet}_{{\cal D}}(X; n)$  is quasiisomorphic to the truncated Deligne complex $\tau_{\leq 2n}\underline \R'_{{\cal D}}(n)(X(\C))$. 
\end{proposition}

{\bf Proof}. See [C] or [G3]. 

Now if $X$ is a variety over $\R$, then
$$
C^{\bullet}_{{\cal D}}(X_{\R}; \R(n)):= C^{\bullet}_{{\cal D}}(X, n)^{\bar F_{{\infty}}}  \qquad  
H_{{\cal D}}^i(X_{\R}; \R(n)) = H^i\Bigl(  C^{\bullet}_{{\cal D}}(X_{\R}; \R(n))\Bigr)
$$

{\bf 2. The Arakelov motivic complexes}. By definition the weight $n$ Arakelov motivic complex 
$\Gamma^{{\cal A}}(X;n)$ is defined as follows
$$
{\rm R}\Gamma^{{\cal A}}_X(n):= 
{\rm Cone}\Bigl({\rm R}\Gamma_X(n) \stackrel{{\rm Reg}}{\lra} {\cal C}_{D}(X_{\R},; \R(n))\Bigr)[-1]
$$
where ${\rm R}\Gamma_X(n)$ is the weight $n$ motivic complex and  
${\rm Reg}$ is the regulator map. 
It is an object of the derived category. When $X$ is regular scheme over $\Q$ then 
$$
{\rm R}\Gamma^{{\cal A}}_{X_{\Q}}(n):= 
{\rm Cone}\Bigl({\rm R}\Gamma_X(n) \stackrel{{\rm Reg}}{\lra} {\cal C}_{D}(X\otimes_{\Q}{\R}; \R(n))\Bigr)[-1]
$$

{\bf Remark}. For the application to special values of $L$-functions we should take a model of $X$ 
over $\Z$ and have a similar construction for finite primes $p$ where it has semistable reduction.

Below we discuss complexes $\Gamma^{\cal A}(X;n)$ with 
the polylogarithmic complex $\Gamma(X;n)$ serving as the motivic complex ${\rm R}\Gamma_X(n)$. 
We will abuse notations by writing 
$({\cal D}^{p,q} \oplus {\cal D}^{q,p})_{\R}(k)$ for the subgroup of $\R(k)$-valued 
currents in ${\cal D}^{p,q} \oplus {\cal D}^{q,p}$. 
The last groups of the complex $\Gamma^{\cal A}(X;n)$ 
look as follows (see also s. 3.6 below):
$$
\begin{array}{cccccc}
... &\oplus_{Y_2 \in X^{(n-2)}}K_2(Y_2) & 
\lra &\oplus_{Y_1 \in X^{(n-1)}}\C(Y_1)^* & \lra &\oplus_{Y \in X^{(n)}}\Z\\
&&&&\\
... &\downarrow r_n(n-2)&&\downarrow r_n(n-1)&&\downarrow r_n(n)\\
&&&&&\\
... & ({\cal D}_{\R}^{n-2,n-1} \oplus  {\cal D}^{n-1,n-2})_{\R}(n-1) &\stackrel{\partial + \overline \partial }{\lra} &{\cal D}^{n-1,n-1}_{\R}(n-1)
&\stackrel{2 \overline \partial \partial}{\lra} &{\cal D}_{cl, \R}^{n,n}(n)
\end{array}
$$
The $(2n)$-th cohomology 
group of this complex looks as follows:
\begin{equation} \label{CHE1}
\frac{\{\mbox{codim. $n$ cycle} \quad Y, \quad g \in {\cal D}^{n-1,n-1}_{cl, \R}(n-1)\} 
\quad  \mbox{such that} \quad 2 \overline \partial \partial g + \delta(Y) =0}
{\{{\rm div}(f), -\log|f|\delta_{Y_1}\}, \quad f \in \C(Y_1)^*; \qquad \{0; \partial u + \overline \partial v\}}
\end{equation}
where $$
u \in {\cal D}_{}^{n-2,n-1}; \quad v \in {\cal D}_{}^{n-1,n-2}
$$
This is the group of  codimension $n$ Arakelov cycles on $X(\C)$. To compare with [S] 
notice that 
$$
d d^{\C} = \frac{1}{2\pi i} \overline \partial \partial, \qquad 2 \overline \partial \partial \log|f| = \delta_{{\rm div} (f)}
$$

If we modify 
this definition by replacing the last group ${\cal D}_{cl, \R}^{n,n}(n)$ by its quotient by smooth forms ${\cal D}_{cl, \R}^{n,n}(n)/{\cal A}_{cl, \R}^{n,n}(n)$  
we get the group of  codimension $n$ Arakelov cycles  defined by Gillet and Soul\'e [S]: 
\begin{equation} \label{CHE2}
\frac{\{Y \in X^{(n)}, \quad g \in {\cal D}_{cl, \R}^{n-1,n-1}(n-1)\} \quad  \mbox{such that} \quad 2 \overline \partial \partial g + \delta(Y)\quad\mbox{is smooth} }
{\{{\rm div}(f), -\log|f|\delta_{Y_1}\}, \quad f \in \C(Y_1)^*; \qquad \{0; \partial u + \overline \partial v\}}
\end{equation}
Such a  modification does not seem  natural in our context, but
it worked well in the Gillet-Soul\'e theory. 
The reason to have several versions of Arakelov groups is the desire to have a theory of 
Chern classes  with values in Arakelov motivic cohomology (and finally the higher 
Riemann-Roch theorem) 
for various versions of the category of vector bundles with some kind of hermitian metric discussed below. 

a) Consider holomorphic 
vector bundles with {\it flat} hermitian metrics. They have Chern classes 
with values in the groups (\ref{CHE1}). In particular  holomorphic line bundles with flat 
hermitian metrics form an abelian group under the tensor product, denoted 
${\widetilde {\rm Pic}}_0(X)$, which sits in an exact sequence 
$$
0 \lra \R^{\pi_0(X(\C))} \lra {\widetilde {\rm Pic}}_0(X) \lra {\rm Pic}_0(X) \lra 0 
$$
Indeed, a flat hermitian hermitian metric in a line bundle over a connected manifold is determined 
up to a constant. 
There is an isomorphism 
$$
{\widetilde {\rm Pic}}_0(X)\quad \stackrel{=}{\lra } \quad H^{2}\Bigl( \Gamma^{{\cal A}}(X;1)\Bigr) 
$$
given by the first Chern class as follows. If $s$ is a section of  a holomorphic 
 line bundle over $X$ with a 
hermitian metric $||\cdot ||$ then the pair 
$(\log||s ||, {\rm div} (s))$ provides the corresponding class in $H^{2}\Bigl( \Gamma^{{\cal A}}(X;1)\Bigr)$. Indeed, according to the Poincar\'e-Lelong formula 
for a meromorphic section $s$ of any holomorphic line bundle one has 
$$
2\overline \partial \partial \log||s || - \delta_{{\rm div} (s)} = c_1(L, ||\cdot||)
$$
where on the right is the Chern form related to the hermitian structure on the line bundle, which is zero if and only if the metric is flat.

b) To have Chern classes for holomorphic vector bundles with 
{\it arbitrary} hermitian metrics we are more or less 
forced to the Gillet-Soul\'e definition. However in this case 
the groups of Arakelov cycles are infinite dimensional. 

c) Choose a Kahler metric on $X(\C)$ and consider only such metrics (harmonic metrics)
 on holomorphic vector bundles over $X(\C)$ whose Chern forms are harmonic with 
respect to the choosen metric. 
Then we have Chern classes yet to another modification of the Arakelov group where ${\cal D}_{cl, \R}^{n,n}(n)$ replaced 
by its quotient by  the image of the following map provided by the Hodge theory:
$$
H^{2n}(X(\C), \Z(n)) \lra {\cal A}_{cl, \R}^{n,n}(n)
$$

{\bf 3. The weight one}. The regulator map on the weight one motivic complex 
looks as follows: 
$$
\begin{array}{ccc}
\C(X)^*& \lra &\oplus_{Y \in X^{(1)}}\Z\\
&&\\
\downarrow r_1(1)&&\downarrow r_1(2)\\
&&\\
{\cal D}^{0,0}&
\stackrel{2 \overline \partial \partial}{\lra} &{\cal D}_{cl, \R}^{1,1}(1)
\end{array}
$$
$$
r_1(2): Y_1 \lms 2 \pi i \cdot \delta_{Y_1}, \qquad 
r_1(1): f \lms \log|f|
$$
Here the top line is the weight $1$ motivic complex, sitting in degrees $[1,2]$. 
Shifting by $-1$ the 
total complex associated with this bicomplex we get the weight $1$ Arakelov 
motivic complex $\Gamma^{\cal A}(X;1)$. 

Here are examples of the regulator maps for the weights $n \leq 3$.  

{\bf 4. The weight two}. The regulator map on the weight two motivic complexes looks as follows. 
$$
\begin{array}{ccccccc}
{\cal B}_2(\C(X))& \stackrel{\delta}{\lra} &\Lambda^2\C(X)^*&\stackrel{\partial}{\lra} &\oplus_{Y \in X^{(1)}}\C(Y)^*&\stackrel{\partial}{\lra}&\oplus_{Y \in X^{(2)}}\Z\\
&&&&&&\\
\downarrow r_2(1)&&\downarrow r_2(2)&&\downarrow r_2(3)&&\downarrow r_2(4)\\
&&&&&&\\
{\cal D}_{\R}^{0,0}(1)& \stackrel{D}{\lra} &({\cal D}^{0,1}\oplus {\cal D}^{1,0})_{\R}(1)&\stackrel{D}{\lra} &{\cal D}_{\R}^{1,1}(1)&
\stackrel{2 \overline \partial \partial}{\lra} &{\cal D}^{2,2}(2)_{cl, \R}
\end{array}
$$
where we set
$$
r_2(4): Y_2 \lms (2 \pi i)^2 \cdot \delta_{Y_2}
$$
$$
r_2(3): (Y_1,f) \lms 2 \pi i \cdot  \log|f|\delta_{Y_1}
$$
$$
r_2(2): f\wedge g \lms -\log|f|d i\arg g + \log|g|d i\arg f 
$$
$$
r_2(1): \{f\}_2 \lms \widehat {\cal L}_2(f)
$$
To prove that we get a morphism of complexes we use theorem \ref{1.11a}. The 
following 
argument is needed to check the commutativity of the second square.
The de Rham differential of the distribution $r_2(2)(f\wedge g)$ is 
$$
{\cal D}  \Bigl(-\log|f|d i\arg g + \log|g|d i\arg f\Bigr) = 
$$
$$
\pi_2(d \log f \wedge \log g) + 2 \pi i \cdot (\log|g| \delta(f) - \log|f| \delta(g)) 
$$ 
This  {\bf does not} coincide with $r_2(3) \circ \partial (f\wedge g)$, but the difference is
$$
({\cal D} \circ r_2(2) - r_2(3) \circ \partial) (f\wedge g) = \pi_2(d \log f \wedge \log g) \in \quad ({\cal D}^{0,2} \oplus {\cal D}^{2,0})_{\R}(1)
$$ 
Defining the differential $D$ on the second group of the complex 
${\cal C}_{{\cal D}}(X_{\R}, \R(2))$ 
we take the de Rham differential  and throw away from it precisely these components. Therefore the middle square is commutative.

{\bf 5. The weight three}. The weight three motivic complex $\Gamma(X;3)$ is the total complex of the following bicomplex:
(the first group is in degree $1$)
$$
\begin{array}{ccccc}
{\cal B}_3(\C(X))& \lra &{\cal B}_2(\C(X)) \otimes \C(X)^*& \lra &\Lambda^3\C(X)^*\\
&&\downarrow &&\downarrow \\
&&\oplus_{Y_1 \in X^{(1)}}{\cal B}_2(\C(Y_1))&\lra &\oplus_{Y_1 \in X^{(1)}}\Lambda^2\C(Y_1)^*\\
&&&&\downarrow \\
&&&&\oplus_{Y_2 \in X^{(2)}}\C(Y_2)^*\\
&&&&\downarrow \\
&&&&\oplus_{Y_3 \in X^{(3)}}\C(Y_3)^*
\end{array}
$$
The Deligne complex 
${\cal C}_{{\cal D}}(X_{\R}, \R(3))$ looks as follows:
$$
\begin{array}{ccccccc}
&&&&&&{\cal D}^{3,3}_{cl, \R}(3)\\
&&&&&2 \overline \partial \partial \nearrow &\\
{\cal D}^{0,2}  & \stackrel{\partial}{\lra} &{\cal D}^{1,2}  &\stackrel{\partial}{\lra} &{\cal D}^{2,2}  &&\\
\uparrow \overline \partial &&\uparrow \overline \partial &&\uparrow \overline \partial &&\\
{\cal D}^{0,1}  &\stackrel{\partial}{\lra} &{\cal D}^{1,1}  &\stackrel{\partial}{\lra} &{\cal D}^{1,2}  &&\\
\uparrow \overline \partial &&\uparrow \overline \partial &&\uparrow \overline \partial &&\\
{\cal D}^{0,0}  &\stackrel{\partial}{\lra} &{\cal D}^{1,0}  &\stackrel{\partial}{\lra} &{\cal D}^{2,0} &&
\end{array}
$$
(Recall that the $3\times 3$ square in this diagram consists of $\R(2)$-valued distributions). 
We construct the regulator map from the motivic complex $\Gamma(X;3)$ to the Deligne complex 
${\cal C}_{{\cal D}}(X_{\R}, \R(3))$
by setting 
$$
r_3(6): Y_3 \lms (2 \pi i)^3 \cdot \delta_{Y_3}
$$
$$
r_3(5): (Y_2,f) \lms (2 \pi i)^2 \cdot  \log|f|\delta_{Y_2}
$$
$$
r_3(4): (Y_1, f\wedge g) \lms 2 \pi i \cdot (-\log|f|d i\arg g + \log|g|d i\arg f )\delta_{Y_1}
$$
$$
r_3(3): (Y_1,\{f\}_2) \lms 2 \pi i \cdot\widehat {\cal L}_2(f)\delta_{Y_1}
$$
$$
r_3(3): f_1\wedge f_2 \wedge f_3 \lms {\rm Alt}_3\Bigl(\frac{1}{6}\log|f_1|d \log|f_2| \wedge d \log|f_3| + \frac{1}{2}\log|f_1|
d i\arg f_2 \wedge    d i\arg f_3\Bigr)
$$
$$
r_3(2): \{f\}_2 \otimes g \lms \widehat {\cal L}_2(f)d i \arg g - \frac{1}{3} \log|g| \cdot 
\Bigl(-\log |1-f| d\log |f| + \log |f| d\log |1-f|\Bigr)
$$
$$
r_3(1): \{f\}_3 \lms \widehat {\cal L}_3(f)
$$

{\bf 6. The general case}. Let $d:= {\rm dim} X$. Then the complex $\Gamma(X;n)$ is the total complex 
of the following bicomplex:
$$
\Gamma(\C(X); n) \lra \oplus_{Y_1 \in X^{(1)}}\Gamma(\C(Y_1); n-1)[-1] \lra 
$$
$$
\oplus_{Y_2 \in X^{(2)}}\Gamma(\C(Y_2); n-2)[-1] \lra ... 
\lra 
\oplus_{Y_d \in X^{(d)}}\Gamma(\C(Y_d); n-d)[-d]
$$
where the arrows are provided by the residue maps, see [G1], p 239-240. 

To define the regulator map from this gadget to the complex ${\cal C}_{{\cal D}}(X_{\R}; \R(n))$ 
we specify it for each of 
$\Gamma(\C(Y_k); n-k)[-k]$ where $k = 0, ..., d$. 
Namely, we take the constructed in theorem 
(\ref{TCON})  homomorphism $r_{n-k}(\cdot)$ for $\eta_{Y_k}$ and multiply it by 
$(2 \pi i)^{n-k} \delta_{Y_{n-k}}$. Notice that the distribution $\delta_{Y}$ 
depends only on the generic point of a subvariety $Y$. 
Its image in the bicomplex described in the section 3.1 is in 
$$
\oplus
{\cal D}^{p, q}; \qquad p \geq k, q \geq k, p+q \leq n+k
$$
which is the triangle  symmetric with respect to the diagonal 
in the $n\times n$ square 
presented on the picture below. 
\begin{center}
\hspace{4.0cm}
\epsffile{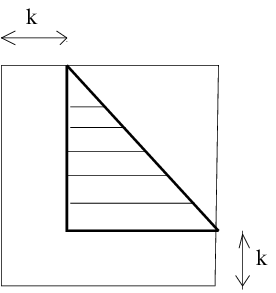}
\end{center}

Compatibility with the differentials is a corollary of theorem 
\ref{TCON} and the following remark. If we would consider the $n\times n$ square in the 
Deligne complex as a part of the Dolbeault complex then we do not get a morphism of complexes. 
For $k<n$ the descrepency lies in the subgroup 
\begin{equation} \label{IIII}
{\cal D}^{k, n+1} \oplus {\cal D}^{n+1, k}
\end{equation}
 and come from the elements 
$(Y_k, f_1\wedge ... \wedge f_{n-k})$. However since the restriction of $D$ to the subgroup 
${\cal D}^{k, n} \oplus {\cal D}^{n, k}$ equals to the de Rham differential ${\cal D}$ 
modulo the components (\ref{IIII}) this does not create a problem.

\section{Proofs}

{\bf 1. The differential equation for the function ${\cal L}_n(z)$}.  
The following proposition was stated without a proof as formula 1.14 in [G2].

\begin{proposition} \label{PPRROO}
The  differential equation for 
$\widehat {\cal L}_{n}(f)$ for $n \geq 3$ is:
\begin{equation} \label{12.4.4}
d \widehat {\cal L}_{n}(f) = \widehat {\cal L}_{n-1}(f) d i \arg z - \sum_{k=2}^{n-1} \beta_k \cdot 
\widehat {\cal L}_{n-k, k}(f)
\end{equation}
\end{proposition}

{\bf Proof}.   Consider the generating series
$$
{Li}(z;t):= \sum_{n \geq 1} {Li}_{n}(z) t^{n-1}; \qquad   \widehat {\cal L}(z;t):= \sum_{n \geq 1} \widehat {\cal L}_{n}(z) t^{n-1}
$$

\begin{lemma} 
Differential equations (\ref{12.4.4}) together with the formulas 
\begin{equation} \label{12.7.1}
d \widehat {\cal L}_1(z) = -d\log|1-z|, \quad d \widehat {\cal L}_2(z) = -\log|1-z|d i \arg z + 
\log|z|d i \arg (1-z) 
\end{equation}
are equivalent to the formula
\begin{equation} \label{12.25.7}
d \widehat {\cal L}(z;t) =   -d\log|1-z| + \log|z| d i \arg (1-z)\cdot t 
+ \widehat {\cal L}(z;t)d  i \arg z \cdot t 
\end{equation}
\begin{equation} \label{12.4.5}
- \left\{ \frac{2 |\log z| \cdot t}{e^{2 \log |z|\cdot t}-1}+ \log|z| \cdot t - 1 \right\} 
\Bigl(\widehat {\cal L}(z;t) \frac{d\log |z|}{\log |z|} + d\log |1-z|\Bigr)\end{equation}
\end{lemma}

{\bf Proof}.  
We check directly that the differential equation (\ref{12.7.1}) appears as the coefficients 
of (\ref{12.4.5}) in $t^0$ and $t^1$. 
From now on we will work modulo the monomials 
$t^0$ and $t^1$. 
The  second term on the right of (\ref{12.4.4}) (written for $n \geq 3$) 
provides the total contribution of
$$
-\sum_{n\geq 3} \Bigl(\sum_{k \geq 2}^{n-2} \beta_k \widehat {\cal L}_{n-k}(z) \cdot \log^{k-1}|z| d \log |z| +
\beta_{n-1} \alpha(1-z,z) \log^{n-2}|z|
\Bigr) t^{n-1} =
$$

$$
-\sum_{2 \leq k < n} \beta_k \widehat {\cal L}_{n-k}(z) \log^{k-1}|z| t^{n-1} \cdot d \log |z| 
$$
$$
 -\sum_{n\geq 3}\beta_{n-1}  \log^{n-1}|z|t^{n-1} \cdot d\log|1-z|\quad = \quad
$$

$$
- \left\{ \frac{2 \log z \cdot t}{e^{2 \log |z|\cdot t}-1} + \log|z| \cdot t - 1\right\} 
\Bigl(\widehat {\cal L}(z;t) \frac{d\log |z|}{\log |z|}  + d\log |1-z| \Bigr)
$$
which is precisely (\ref{12.4.5}). The  first term on the right of (\ref{12.4.4}) gives the last term in (\ref{12.25.7}). 
The lemma is proved. 

Let us prove formula (\ref{12.25.7}) - (\ref{12.4.5}). 
One has for $n \geq 2$:
$$
d \pi_n(Li_m(z)) =  \frac{1}{2}d\Bigl(Li_m(z) +(-1)^{n-1}\overline{Li_m(z)}\Bigr ) =
$$
\begin{equation} \label{12.3.1}
\pi_n (Li_{m-1}(z)) d \log|z| + \pi_{n-1} (Li_{m-1}(z)) di \arg z
\end{equation}

Define an operator $P$ 
acting on the generating series for a sequence of holomorphic functions $f_n(z)$ by 
$$
P\Bigl( \sum_{n\geq 1}f_n(z)t^{n-1}\Bigr):= \sum_{n\geq 1} \pi_n(f_n) \cdot t^{n-1}
$$
Then 
$$
\widehat {\cal L}(z;t) = P\left\{\sum_{0 \leq k < n}
\beta_k Li_{n-k}(z) \cdot \log^k|z| \cdot t^{n-1}\right\} = 
$$
\begin{equation} \label{12.4.1}
P\left\{Li(z;t) \cdot \frac{2 \log |z| t}{e^{2 \log |z| t}-1} \right\} 
\end{equation}

One has 
$$
d\Bigl( \frac{2x}{e^{2x}-1}\Bigr) = \Bigl( \frac{2x}{e^{2x}-1}    
- \frac{4x^2e^{2x}}{(e^{2x}-1)^2}\Bigr) \frac{dx}{x} = \frac{2x}{e^{2x}-1} \cdot \Bigl( 1-2x -\frac{2x}{e^{2x}-1}\Bigr)\frac{dx}{x}
$$
Thus since the second factor in (\ref{12.4.1}) is a real function we get 
$$
d \widehat {\cal L}(z;t) = d_1 P\Bigl(Li(z;t)
\frac{2\log |z| \cdot t }{e^{2 \log |z|\cdot t}-1} \Bigr)  
$$
$$
+P\left\{{Li}(z;t) \cdot 
\frac{2\log |z| \cdot t }{e^{2 \log |z|\cdot t}-1}  \cdot 
\Bigl(1 -   2\log |z| \cdot t -   \frac{2\log |z| \cdot t }{e^{2 \log |z|\cdot t}-1}\Bigr)\right\} 
\frac{d\log |z| }{\log |z|} 
$$
where $d_1$ is the differential applied to  $Li(z;t)$ only. Since 
$\frac{2x}{e^{2x}-1}+x-1$ is an even function in $x$ we rewrite this as 
\begin{equation} \label{12.25.1}
d_1 P\Bigl(Li(z;t)\frac{2\log |z| \cdot t }{e^{2 \log |z|\cdot t}-1} \Bigr)   
\end{equation}
\begin{equation} \label{12.25.8}
- \widehat {\cal L}(z;t)\cdot 
\Bigl(\frac{2\log |z| \cdot t }{e^{2 \log |z|\cdot t} -1} + \log |z|\cdot t -1\Bigr)\frac{d\log |z| }{\log |z|} 
\end{equation}
\begin{equation} \label{12.25.2}
-P\left\{ Li(z;t)\frac{2\log |z| \cdot t }{e^{2 \log |z|\cdot t}-1} \cdot t\right\}d\log |z|
\end{equation}
We handle this as follows. 

i) (\ref{12.25.8}) matches half of (\ref{12.25.7}). 

ii) The $Li_1$-part of  (\ref{12.25.1}) is 
$$
d_1 P\Bigl(Li_1(z)\frac{2\log |z| \cdot t }{e^{2 \log |z|\cdot t}-1} \Bigr)  \quad = \quad - \frac{2\log |z| \cdot t }{e^{2 \log |z|\cdot t}-1} d \log |1-z|
$$
It matches the  second half of  (\ref{12.4.5}) modulo coefficients in $t^0$ and $t^1$. 

iii) The rest of (\ref{12.25.1}) equals to 
\begin{equation} \label{12.25.3}
P\left\{ Li(z;t)\frac{2\log |z| \cdot t }{e^{2 \log |z|\cdot t}-1} d\log z \cdot t\right\} = 
\end{equation}
\begin{equation} \label{12.25.4}
P\left\{ \widehat {\cal L}(z;t)(z;t)\right\}d i\arg z \cdot t 
\quad + \quad P\left\{ Li(z;t)\frac{2\log |z| \cdot t }{e^{2 \log |z|\cdot t}-1} \cdot t\right\} d\log |z|
\end{equation}
Thus the first term in (\ref{12.25.4}) matches the right term in (\ref{12.25.7}), and 
the second is canceled with (\ref{12.25.2}). The proposition is proved.

{\bf 2.       Proof of theorem \ref{TCON}}. i) Let us show that the map
\begin{equation} \label{1.15}
\{f\} \lms \widehat {\cal L}_{n-k,k}(f) 
\end{equation}
provides a group homomorphism 
$
{\cal B}_n(\C(X)) \lra {\cal A}_{\eta}^1(X(\C))
$. 
The arguments are similar to the proof of theorem 1.15 in [G2]. Suppose first that $n-k>1$. 
By theorem  1.15 in [G2] the map 
$
\{f\} \lms \widehat {\cal L}_{n}(f) 
$ 
provides a group homomorphism ${\cal B}_n(\C(X)) \lra {\cal A}_{\eta}^0(X(\C))$. 
Consider the following map 
\begin{equation} \label{1.16}
{\cal B}_n(\C(X)) \lra {\cal B}_{n-k}(\C(X)) \otimes S^k\C(X)^*, \quad \{f\}_n \lms \{f\}_{n-k} \otimes f^k 
\end{equation}
It  can be defined as a composition 
$$
\{f\}_n \lms \{f\}_{n-1} \otimes f  
\lms \{f\}_{n-2} \otimes f\cdot f   \lms ... \lms \{f\}_{n-k} \otimes f\cdot ... \cdot f 
$$
Each of these maps is of type $\delta \otimes id$ and thus is a homomorphism of abelian groups. 
The composition followed by the homomorphism 
$$
{\cal B}_{n-k}(\C(X)) \otimes S^k\C(X)^* \lra {\cal A}_{\eta}^1(X(\C))
$$
$$
\{f\}_{n-k} \otimes g_1 \cdot ... \cdot g_k \lms \widehat {\cal L}_{n-k}(f) \frac{1}{k}d(\log |g_1| \cdot ... \cdot\log |g_k|)  
$$
leads to the map (\ref{1.15}).  In the case $n-k=1$ we present map (\ref{1.15}) as a composition 
$$
{\cal B}_n(\C(X)) \stackrel{(\ref{1.16})}{\lra} {\cal B}_{2}(\C(X)) \otimes S^{n-2}\C(X)^* 
\stackrel{a}{\lra}{\cal A}_{\eta}^1(X(\C))
$$
$$
a: \{f\}_2 \otimes g_1 \cdot ... \cdot g_{n-2} \lms \alpha(1-f, f) \cdot \log|g_1| \cdot ... \cdot \log|g_{n-2}|
$$
and use the fact that $\{f\}_2 \lms \alpha(1-f, f)$ is a group homomorphism. 

ii) The properties (a) and (b) in theorem \ref{CANHH} are 
true by the very definitions. 

iii) It is very easy to see that the property (d) is also true by the very definitions: 
it basically reflects the fact 
that the numerical coefficients in the formula for $r_{n+m}(m+1)$ do not depend on $m$ 
and uses formula (\ref{3/5/00.6}).  

iv) It remains to prove that the diagram in theorem \ref{CANHH} is commutative. Let us show first that its right square 
$$
\begin{array}{ccc}
{\cal B}_2(\C(X)) \otimes \Lambda^{n-2}\C(X)^* & \stackrel{\delta}{\lra} & \Lambda^n\C(X)^*\\
&&\\
\downarrow r_{n}(n-1)& &\downarrow r_n(n)\\
&&\\
{\cal A}^{n-2}_{\eta}(X(\C)) & \stackrel{d}{\lra} & {\cal A}^{n-1}_{\eta}(X(\C))
\end{array}
$$
 is commutative. We change the numeration by putting $n:= m+2$.

\begin{proposition} \label{n=2,m}
$$
r_{m+2}(m+1): \{f\}_2 \otimes g_1 \wedge ... \wedge g_m \lms 
$$
$$
\widehat {\cal L}_{2}(f) \cdot {\cal A}_m 
\left\{  \sum_{p \geq 0 } \frac{1}{2p+1} \bigwedge_{i=1}^{2p}d\log|g_i| \wedge 
\bigwedge_{j=2p+1}^{m}d i \arg g_{j}  \right\}
$$
$$
- \alpha(1-f, f) \wedge {\cal A}_m \Bigl(
\sum_{p \geq 0 } \frac{1}{(2p+1)(2p+3)} \log|g_1| \cdot \bigwedge_{i=2}^{2p+1}d\log|g_i| \wedge 
\bigwedge_{j=2p+2}^{m}d i \arg g_{j}  \Bigr)
$$
\end{proposition}

{\bf Proof}. It follows immediately from the definition and lemma \ref{RADEMA}:

Let us compute
$$
d \Bigl( r_{m+2}(m+1) (\{f\}_2 \otimes g_1 \wedge ... \wedge g_m) \Bigr) 
$$
Using proposition \ref{n=2,m} we get
\begin{equation} \label{TG1}
\Bigl( - \log |1-f| d i \arg f + \log |f| d i \arg (1-f) \Bigr) \wedge 
\end{equation}
$$
{\cal A}_m 
\left\{  \sum_{p \geq 0 } \frac{1}{2p+1} \bigwedge_{i=1}^{2p}d\log|g_i| \wedge 
\bigwedge_{j=2p+1}^{m}d i \arg g_{j} 
  \right\}
$$

$$
+ \quad 2 d\log |1-f|\wedge d\log|f| \wedge 
$$
\begin{equation} \label{TG2}
{\cal A}_m \Bigl(
\sum_{p \geq 0 } \frac{1}{(2p+1)(2p+3)} \log|g_1|  \cdot \bigwedge_{i=2}^{2p+1}d\log|g_i| \wedge 
\bigwedge_{j=2p+2}^{m}d i \arg g_{j}\Bigr)
\end{equation}

\begin{equation} \label{TG3}
+ \quad \alpha(1-f,f) \wedge {\cal A}_m \Bigl(
\sum_{p \geq 0 } \frac{1}{2p+3} \bigwedge_{i=1}^{2p+1}d\log|g_i| \wedge 
\bigwedge_{j=2p+2}^{m}d i \arg g_{j}\Bigr)
\end{equation}

We can write the differential form $r_m(g_1\wedge ...\wedge g_m)$ as follows:
\begin{equation} \label{TG4}
-{\cal A}_m\left\{\sum_{p \geq 0}\frac{1}{2p+1}\log|g_1| \cdot 
\bigwedge_{i=2}^{2p+1}d\log|g_i| \wedge 
\bigwedge_{j=2p+2}^{m}d i \arg g_{j}   \right\}
\end{equation}
In particular this makes transparent the property
$$
d r_m(g_1 \wedge  ...\wedge  g_m) + \pi_m(d \log g_1 \wedge ... \wedge  \log g_m) = 0
$$

Now computing
\begin{equation} \label{TG5}
d r_{m+2}((1-f) \wedge  f \wedge g_1 \wedge  ... \wedge  g_m)
\end{equation}
using formula (\ref{TG4}) and comparing the result with formulas (\ref{TG1}) - (\ref{TG3})  we see the following:

1) Formula (\ref{TG1}) matches the part of (\ref{TG5}) where the contribution of $(1-f)\wedge f$ is given by 
$-\log|1-f| d i \arg|f| + \log|f| d i\arg|1-f|$.

2) Formula (\ref{TG3}) matches the part of (\ref{TG5}) where the contribution of $(1-f)\wedge f$ is given by 
$-\log|1-f| d \log |f| + \log|f| d \log |1-f|$.

Before we continue any further let us note that 
\begin{equation} \label{TG7}
d \log|1-f| \wedge d i\arg (f) \quad = \quad d\log|f|  \wedge di\arg (1-f) 
\end{equation}
\begin{equation} \label{TG8}
d \log|1-f| \wedge d \log|f| \quad = -\quad di\arg (1-f)  \wedge d i\arg f 
\end{equation}
Indeed,  
$$
 0 = d \log (1-f) \wedge d \log f  = \Bigl(d \log |1-f| + d i \arg (1-f) \Bigr)\wedge \Bigl(d \log |f| + d \arg i f\Bigr)
$$ 

Therefore using (\ref{TG7}) we see the following: 

3) the sum of the terms of (\ref{TG5}) where the contribution of $(1-f)\wedge f$ is given either by 
$d\log|1-f| \wedge d i\arg |f|$ or by $ d i \arg (1-f)f\wedge  d \log |f|$ is zero. 

4) It remains to show that (\ref{TG2}) matches the part of (\ref{TG5}) where 
the contribution of $(1-f)\wedge f$ is given by 
$d\log|1-f| \wedge d \log |f|$ or $ d i \arg (1-f)f\wedge  d i\arg f$.

 Let us substitute
$$
\frac{2}{(2p+1)(2p+3)} = \frac{1}{2p+1} - \frac{1}{2p+3} 
$$
into (\ref{TG2}) and split the formula we get into two parts. The first part, denoted $(\ref{TG2})_I$, is the one where 
the terms appear with the coefficient $\frac{1}{2p+1}$. The second part, denoted $(\ref{TG2})_{II}$,  is the rest. 
It corresponds to $\frac{1}{2p+3}$. 

Using (\ref{TG8}) to calculate the part of (\ref{TG5}) where 
$(1-f)\wedge f$ contributes via $ d i \arg (1-f)f\wedge  d i\arg f$ we 
see that it matches with $(\ref{TG2})_I$.  
The other part, $(\ref{TG2})_{II}$,  matches the part of (\ref{TG5}) where 
$(1-f)\wedge f$ contributes $ d \log | 1-f|\wedge  d\log  |f|$. So we 
achieved our  goal 4). 
The commutativity of the last square is proved. 

Now let us prove the commutativity of the squares different from the last one, i.e. ($n>2$)
$$
\begin{array}{ccc}
{\cal B}_n(\C(X)) \otimes \Lambda^{m}\C(X)^* & \stackrel{\delta}{\lra} & 
{\cal B}_{n-1}(\C(X)) \otimes \Lambda^{m+1}\C(X)^*\\
&&\\
\downarrow r_{n+m}(m+1)& &\downarrow r_{n + m}(m+2)\\
&&\\
{\cal A}^{m}_{\eta}(X(\C)) & \stackrel{d}{\lra} & {\cal A}^{m+1}_{\eta}(X(\C))
\end{array}
$$
One has 

\begin{equation} \label{12.26.40}
d r_{n+m}(m+1)\Bigl(\{f\}_n \otimes g_1 \wedge ... \wedge g_m \Bigl)= 
\end{equation}
\begin{equation} \label{12.26.3}
d \widehat {\cal L}_{n}(f) \wedge {\cal A}_m\left\{  \sum_{p \geq 0 } \frac{1}{2p+1} 
\bigwedge_{i=1}^{2p}d\log|g_i| \wedge   \bigwedge_{j=2p+1}^{m}d i \arg g_{j}\right\} 
\end{equation}

\begin{equation} \label{12.26.2}
+ \sum_{k \geq 1 }\sum_{1 \leq p \leq m} \beta_{k, p} d \widehat {\cal L}_{n-k,k}(f) \wedge 
{\cal A}_m\left\{\log|g_1| \bigwedge_{i=2}^{p}d\log|g_i| \wedge   
\bigwedge_{j=p+1}^{m}d i \arg g_{j}
\right\}\end{equation} 

\begin{equation} \label{12.26.1}
-\sum_{k \geq 1 }\sum_{1 \leq p \leq m} p \cdot \beta_{k, p} \widehat {\cal L}_{n-k,k}(f) \wedge 
{\cal A}_m \left\{\bigwedge_{i=1}^{p}d\log|g_i| \wedge   \bigwedge_{j=p+1}^{m}d i \arg g_{j}
 \right\}
\end{equation}

Let us compare this expression with 
\begin{equation} \label{12.26.30}
r_{n+m}(m+2)\Bigl(\{f\}_{n-1} \otimes f \wedge g_1 \wedge ... \wedge g_m \Bigl) = (\ref{12.26.11})_{m+2} + (\ref{12.26.12})_{m+2}
\end{equation}
where $(\ref{12.26.11})_{m+2}$ (resp. $(\ref{12.26.12})_{m+2}$) stays for the term similar
 to  $(\ref{12.26.11})$ (resp. $(\ref{12.26.12})$) in the definition of  $r_{n+m}(m+1)$ (see 
theorem-construction \ref{TCON}). 

Observe that {\it a priori} 
$f$ can contribute to $(\ref{12.26.11})_{m+2}$ via $d\log|f|$ or $d i \arg f$, and to 
$(\ref{12.26.12})_{m+2}$via $\log|f|$, $d\log|f|$ or $d i \arg f$. A careful reader should
 keep track of these five different cases during the proof. 

We split the job in two parts: 

A) Check that the diagram is commutative if we pay attention {\it only} to the terms where $d\log|1-f|$ or $di\arg (1-f)$ appears. We call such terms $A$-terms.

B) Then we check that the diagram is commutative modulo the 
terms with $d\log|1-f|$ or $di\arg (1-f)$.

{\bf Part A)}.  
Here is the contribution from (\ref{12.26.3}) - (\ref{12.26.1}). Recall
that $n \geq 3$. 

From (\ref{12.26.3}) using (\ref{12.4.4}) we get:
\begin{equation} \label{12.99.1}
\beta_{n-1} \cdot \log^{n-1}|f| d\log|1-f| \wedge 
{\cal A}_m\left\{  \sum_{p \geq 0 } \frac{1}{2p+1} 
\bigwedge_{i=1}^{2p}d\log|g_i| \wedge   \bigwedge_{j=2p+1}^{m}d i \arg g_{j}\right\} 
\end{equation}

From (\ref{12.26.2}): 
\begin{equation} \label{12.99.3}
 \sum_{k = 1 }^{n-3} \sum_{1 \leq p \leq m} - \beta_{k, p} \beta_{n-k-1} \log^{n-2}|f| d\log|1-f| 
\wedge d\log|f| \wedge 
\end{equation} 
\begin{equation} \label{12.99.4}
{\cal A}_m\left\{\log|g_1| \bigwedge_{i=2}^{p}d\log|g_i| \wedge   
\bigwedge_{j=p+1}^{m}d i \arg g_{j} 
\right\}
\end{equation} 

\begin{equation} \label{12.99.5}
+\sum_{1 \leq p \leq m} \beta_{n-2, p} d i \arg (1-f)\wedge d\log|f| \cdot 
\log^{n-2}|f| \wedge (\ref{12.99.4}) 
\end{equation}
\begin{equation} \label{12.99.6}
-\sum_{1 \leq p \leq m} (n-1)\beta_{n-1, p} d\log|1-f| 
\wedge d\log|f| \wedge \log^{n-2}|f|\wedge (\ref{12.99.4}) 
\end{equation}

From (\ref{12.26.1}):
\begin{equation} \label{12.99.7}
\sum_{1 \leq p \leq m} p\beta_{n-1, p}  
\log^{n-1}|f| \cdot d\log|1-f| \wedge 
\end{equation}
\begin{equation} \label{12.99.8}
{\cal A}_{m} \left\{\bigwedge_{i=1}^{p}d\log|g_i| \wedge   \bigwedge_{j=p+1}^{m}d i \arg g_{j}
 \right\}
\end{equation}

The contribution of $(\ref{12.26.11})_{m+2}$ is zero. 
The factor $f$ from $\{f\}_{n-1} \otimes f \wedge g_1 \wedge ... g_m $ can contribute into $(\ref{12.26.12})_{m+2}$ via  $\log |f|$, $d i\arg f$ and $d \log |f|$. 

The terms where $f$ contributes via $\log |f|$ are 
\begin{equation} \label{12.101.8}
\sum_{k \geq 1}\sum_{1 \leq p \leq m} \beta_{k, p+1} \widehat {\cal L}_{n-1-k,k}(f) \log|f| 
\wedge 
(\ref{12.99.8})
\end{equation}

Since we are concerned only with the $A$-terms, this boils down to
\begin{equation} \label{12.100.8}
\sum_{1 \leq p \leq m} \beta_{n-2, p+1} \widehat {\cal L}_{1, n-2}(f) \log|f| \wedge (\ref{12.99.8}) = 
\end{equation}
$$
-\sum_{1 \leq p \leq m} \beta_{n-2, p+1}  \log^{n-1}|f| d\log|1-f|\wedge (\ref{12.99.8})
$$
Using formulas (\ref{12.99}) we see that it matches (\ref{12.99.1}) + (\ref{12.99.7}). 

 The terms where $f$ contributes into $(\ref{12.26.12})_{m+2}$ via $d i\arg f$ are 
\begin{equation} \label{12.102.8}
-\sum_{k \geq 1}\sum_{1 \leq p \leq m} \beta_{k, p} \widehat {\cal L}_{n-1-k,k}(f) d i \arg f 
\wedge 
(\ref{12.99.4})
\end{equation}
Its contribution modulo $B$-terms is 
\begin{equation} \label{12.102.8}
-\sum_{1 \leq p \leq m} \beta_{n-2, p} \log^{n-2}|f| d\log|1-f|\wedge  d i \arg f 
\wedge 
(\ref{12.99.4})
\end{equation}
Using (\ref{TG7}) we see that it matches (\ref{12.99.5}). 

Similar considerations for terms where $f$ contributes into $(\ref{12.26.12})_{m+2}$ via 
$d \log| f|$ minus ((\ref{12.99.6}) + (\ref{12.99.3})) leads  to the left hand side of 
(\ref{3/5/00.1}), which is zero according to proposition \ref{3/5/00.3}.  
The part A) is proved.

{\bf Part B)}. 
We will write $ X\stackrel{A}{=} Y$ if $X-Y = 0 $ modulo the $A$-terms.

a) Notice that the contribution of $f$ to $(\ref{12.26.12})_{m+2}$ via $d\log|f|$ 
is zero since $\widehat {\cal L}_{n-k,k}(f) \wedge  d\log|f| = 0$ modulo $A$-terms.

b) The contribution to (\ref{12.26.3}) of the term $\widehat {\cal L}_{n-1}(f) d i\arg f$ 
in  formula (\ref{12.4.4}) for $d \widehat {\cal L}_{n}(f) $ 
matches  the part of term $(\ref{12.26.11})_{m+2}$ where $f$ contributes via $di\arg f$. 

c) One has
\begin{equation} \label{12.27.1}
d\widehat {\cal L}_{n-k,k}(f) \quad \stackrel{A}{=} \quad - \widehat {\cal L}_{n-k-1,k}(f) \wedge d i \arg f 
\end{equation}
Thus (\ref{12.26.2}) matches the part of $(\ref{12.26.12})_{m+2}$ where $f$ contributes via 
$d i \arg f$.

d) The contribution to (\ref{12.26.3}) of the other term in  $d \widehat {\cal L}_{n}(f) $ 
 modulo $A$ is
\begin{equation} \label{12.26.20}
-\sum_{p \geq 0 }\sum_{k=2 }^{n-2} \frac{\beta_k}{2p+1} \widehat {\cal L}_{n-k,k-1}(f) \wedge 
\end{equation}
\begin{equation} \label{12.26.21}
{\cal A}_m\left\{   
\log |f| \cdot \bigwedge_{i=1}^{2p}d\log|g_i| \wedge   \bigwedge_{j=2p+1}^{m}d i \arg g_{j}\right\} 
\end{equation}

The only other terms in (\ref{12.26.40})  
with factor (\ref{12.26.21}) are the terms of 
(\ref{12.26.1}) where $k>1$ and the number of $d\log |g_i|$'s is even, 
i.e.  
\begin{equation} \label{12.26.32}
-\sum_{k \geq 1}2p \beta_{k,2p} \widehat {\cal L}_{n-k,k-1}(f) \wedge (\ref{12.26.21})
\end{equation}
On the other hand, in (\ref{12.26.30}) the terms with factor (\ref{12.26.21}) are
\begin{equation} \label{12.26.33}
\sum_{k \geq 1} \beta_{k,2p+1} \widehat {\cal L}_{n-1-k,k}(f) \wedge (\ref{12.26.21})
\end{equation}
which are precisely the terms of $(\ref{12.26.12})_{m+2} $ where $f$ contributes via $\log|f|$, 
and $p$ in 
$\beta_{k,p} $ is odd. 

Thus out of (\ref{12.26.20}), (\ref{12.26.32}) and (\ref{12.26.33}) we see that 
$\widehat {\cal L}_{n-k,k-1}(f) \wedge(\ref{12.26.21})$ appears in our check of 
 the commutativity of the diagram with factor 
$$
\beta_{k-1,2p+1} + \frac{\beta_{k}}{2p+1}  +2p \beta_{k,2p} \quad \stackrel{(\ref{12.99})}{=} \quad 0
$$

e) The term of $(\ref{12.26.11})_{m+2} $ where $f$ contributes via $d\log|f|$, i.e.
\begin{equation} \label{1999.1}
\sum_{p \geq 0 }\frac{1}{2p+1} \widehat {\cal L}_{n-1}(f) \cdot {\cal A}_{m+1}\left\{   
d\log|f| \wedge \bigwedge_{i=1}^{2p-1}d\log|g_i| \wedge   \bigwedge_{j=2p}^{m}d i \arg g_{j}\right\} 
\end{equation}
matches the term (\ref{12.26.1}) with $k=1$, $p$: odd, i.e.
\begin{equation} \label{1999.2}
-\sum_{p \geq 0 }(2p-1) \beta_{1,2p-1}   \widehat {\cal L}_{n-1,1}(f) \wedge {\cal A}_{m}\left\{   
\bigwedge_{i=1}^{2p-1}d\log|g_i| \wedge   \bigwedge_{j=2p}^{m}d i \arg g_{j}\right\} 
\end{equation}
Indeed, according to (\ref{1999.5}) $(2p-1) \beta_{1,2p-1} = -\frac{1}{2p+1}$. 

f) All the terms of (\ref{12.26.1}) with $k=1$ and even $p$  are zero since $\beta_{1,2p}=0$ by 
(\ref{1999.5}).  

g) The terms of (\ref{12.26.1}) with odd $p$ and $k>1$ look as follows:
\begin{equation} \label{1999.3}
\sum_{p \geq 0 }\sum_{k \geq 2 } (2p+1) \beta_{k,2p+1}   \widehat {\cal L}_{n-k,k-1}(f) 
\wedge 
\end{equation}
\begin{equation} \label{1999.4`}
{\cal A}_{m+1}\left\{  \log|f|  
\bigwedge_{i=1}^{2p+1}d\log|g_i| \wedge   \bigwedge_{j=2p+2}^{m}d i \arg g_{j}\right\} 
\end{equation}
They match with the terms of $(\ref{12.26.12})_{m+2} $  where $f$ contributes via $\log|f|$ and $p$ is even, see (\ref{12.101.8}). 

The  the part B) of the commutativity proof is now complete. 
To double check that 
we list all the possibilities for contribution of $f$ in (\ref{12.26.30}),  
and indicate for each of them the part of the proof where it was handled. 

Contribution of $f$ to $(\ref{12.26.11})_{m+1} $: 

via $d\log|f|$: see e);$\qquad$
via $di\arg f$: see b).

Contribution of $f$ to $(\ref{12.26.12})_{m+2} $: 

via $\log|f|$, $p$ odd: see d); $\qquad$
via $\log|f|$, $p$ even: see g);

via $d\log|f|$: see a); $\qquad$
via $di\arg f$: see c).

Contribution of $f$ to (\ref{12.26.3}): see b) and d).

Contribution of $f$ to (\ref{12.26.2}): see c).

Contribution of $f$ to (\ref{12.26.1}): 

$k=1$, $p$ odd: see e);    $\qquad$        $k>1$, $p$ odd: see g);

$k=1$, $p$ even: see f);   $\qquad$         $k>1$, $p$ even: see d).

The main theorem is proved.


\bigskip

Department of mathematics, Brown University, Providence, RI 02912.
e-mail sasha@math.brown.edu

\end{document}